\documentclass[11pt]{article}
\def\thetitle{A very robust Ramsey theorem for hypergraph matchings}

\usepackage{amsmath,amssymb}
\usepackage{needspace}

\usepackage[usenames,dvipsnames,svgnames,table]{xcolor}
\definecolor{CombinatoricaAqua}{HTML}{00698C}
\definecolor{CombinatoricaBlue}{HTML}{3A3293}
\definecolor{CombinatoricaBrown}{HTML}{66220C}
\definecolor{HarvardCrimson}{rgb}{0.6471,0.1098,0.1882}

\makeatletter
\let\reftagform@=\tagform@
\def\tagform@#1{\maketag@@@{(\ignorespaces\textcolor{CombinatoricaBrown}{#1}\unskip\@@italiccorr)}}
\renewcommand{\eqref}[1]{\textup{\reftagform@{\ref{#1}}}}
\makeatother

\usepackage[backref=page]{hyperref}
\hypersetup{
  unicode,
  pdfencoding=auto,
  pdfauthor={Peter Keevash and Peleg Michaeli},
  pdftitle={\thetitle},
  colorlinks=true,
  citecolor=CombinatoricaBlue,
  linkcolor=CombinatoricaAqua,
  anchorcolor=CombinatoricaBrown,
  urlcolor=HarvardCrimson
}

\usepackage{amsthm}
\usepackage{thmtools}
\makeatletter
\IfFormatAtLeastTF{2024-11-01}{
  \renewcommand*\@addtoreset[2]{%
    \bgroup
      \edef\aliasctr@@truelist{\aliasctr@follow{#2}}%
      \let\@elt\relax
      \expandafter\@cons\aliasctr@@truelist{{#1}}%
    \egroup
    \expandafter\xdef\csname theH#1\endcsname{%
      \expandafter\noexpand\csname theH#2\endcsname.%
      \noexpand\the\noexpand\value{#1}}%
  }
}{}
\makeatother
\usepackage[capitalize]{cleveref}
\declaretheoremstyle[
  spaceabove=\topsep,
  spacebelow=\topsep,
  headfont=\color{CombinatoricaBrown}\normalfont\bfseries,
  bodyfont=\itshape
]{thm}
\declaretheorem[style=thm,parent=section]{theorem}
\declaretheorem[style=thm,sibling=theorem]{lemma}
\declaretheorem[style=thm,sibling=theorem]{proposition}

\declaretheorem[style=thm,sibling=theorem]{conjecture}

\usepackage[nobysame,msc-links,non-sorted-cites]{amsrefs}

\BibSpec{book}{%
  +{}  {\PrintPrimary}                {transition}
  +{,} { \textbf}                     {title}
  +{.} { }                            {part}
  +{:} { \textit}                     {subtitle}
  +{,} { \PrintEdition}               {edition}
  +{}  { \PrintEditorsB}              {editor}
  +{,} { \PrintTranslatorsC}          {translator}
  +{,} { \PrintContributions}         {contribution}
  +{,} { }                            {series}
  +{,} { \voltext}                    {volume}
  +{,} { }                            {publisher}
  +{,} { }                            {organization}
  +{,} { }                            {address}
  +{,} { \PrintDateB}                 {date}
  +{,} { }                            {status}
  +{}  { \parenthesize}               {language}
  +{}  { \PrintTranslation}           {translation}
  +{;} { \PrintReprint}               {reprint}
  +{.} { }                            {note}
  +{.} {}                             {transition}
  +{}  {\SentenceSpace \PrintReviews} {review}
}
\makeatletter
\renewcommand{\PrintNames@a}[4]{%
  \PrintSeries{\name}
  {#1}
  {}{ and \set@othername}
  {,}{ \set@othername}
  {}{ and \set@othername}
  {#2}{#4}{#3}%
}
\makeatother

\usepackage{sectsty}
\sectionfont{\color{CombinatoricaBrown}}
\subsectionfont{\color{CombinatoricaBrown}}

\usepackage[a4paper]{geometry}
\usepackage{mleftright}
\mleftright
\newcommand{\cH}{\mathcal{H}}
\newcommand{\bG}{\mathbf{G}}
\newcommand{\bU}{\mathbf{U}}
\newcommand{\vbeta}{\boldsymbol{\beta}}
\newcommand{\pr}{\mathbb{P}}
\newcommand{\E}{\mathbb{E}}
\newcommand{\eps}{\varepsilon}
\newcommand{\NN}{\mathbb{N}}
\newcommand{\vt}{\mathbf{t}}
\newcommand{\vk}{\mathbf{k}}
\newcommand{\tmax}{\|\vt\|_\infty}
\newcommand{\afl}{\mathsf{AFL}}
\newcommand{\one}{\mathbf{1}}
\newcommand{\defn}[1]{{\bfseries #1}}
\newcommand{\whp}{\textbf{whp}}
\newcommand{\Whp}{\textbf{Whp}}
\newcommand{\sm}{\smallsetminus}
\newcommand{\floor}[1]{\left\lfloor{#1}\right\rfloor}
\newcommand{\ceil}[1]{\left\lceil{#1}\right\rceil}
\DeclareMathOperator{\diam}{diam}

\title{\thetitle}
\author{
  Peter Keevash\thanks{
    Mathematical Institute,
    University of Oxford,
    Oxford, UK.
    Email: \href{mailto:keevash@maths.ox.ac.uk}
                {\tt keevash@maths.ox.ac.uk}.
    Supported by ERC Advanced Grant 883810.
  }
  \and
  Peleg Michaeli\thanks{
    Mathematical Institute,
    University of Oxford,
    Oxford, UK.
    Email: \href{mailto:michaeli@maths.ox.ac.uk}
                {\tt michaeli@maths.ox.ac.uk}.
    Research supported by UKRI under the Horizon Europe Guarantee for MSCA Postdoctoral Fellowships (EP/Z001781/1, IRIS).
  }
}
\date{}

\begin{document}

\maketitle
\thispagestyle{empty}

\begin{abstract}
A fundamental result in Ramsey theory due to Alon, Frankl and Lov\'asz gives an exact formula
for the $q$-colour Ramsey numbers of matchings in $r$-uniform hypergraphs:
if $n \ge \afl_r(\vt)=(r-1)\tmax+1+\sum_{j=1}^q(t_j-1)$ then in any $q$-edge-colouring of $K^r_n$
there is some colour $j$ with a matching of size $t_j$. We show that a very weak
quasirandomness assumption suffices for this Ramsey property: we can replace $K^r_n$ 
by any $r$-graph $H$ on at least $ \afl_r(\vt)+o(n)$ vertices that is a sublinear connector,
meaning that it has no balanced $r$-partite holes of linear size.
This generalises the result for $r=2$ (ordinary graphs) that we proved previously
using a compression algorithm based on Gallai-Edmonds decompositions.
This structure theory is unavailable in the hypergraph setting, so we combine two other proof methods, namely
(a) simultaneous optimisation in all colours by clustering vertices using fractional cover weights from LP duality,
and (b) local-to-global structure via random sampling and Kahn's matching theorem.
\end{abstract}

\section{Introduction}

Ramsey theory studies the monochromatic configurations that must occur in colourings of sufficiently large structures.
In particular, Ramsey's theorem~\cite{Ramsey30} guarantees that every edge-colouring
of a sufficiently large complete uniform hypergraph with a given number of colours
contains a monochromatic copy of any prescribed finite hypergraph of the same uniformity.
Determining the minimum number of vertices needed is a central problem in Ramsey theory.
We will be concerned with the multicolour Ramsey numbers of matchings,
which are known exactly, even when different matching sizes are prescribed in different colours.

We use the following notation.
Given $r$-uniform hypergraphs (\defn{$r$-graphs} for short)
$F_1,\dots,F_q$, we write $H\to(F_1,\dots,F_q)$
if every $q$-edge-colouring of $H$ has a $j$-coloured copy of $F_j$ for some $j\in[q]:=\{1,\dots,q\}$.
Write $\NN$ for the set of non-negative integers
and $\NN_+ = \NN \sm \{0\}$.
For $\vt=(t_1,\dots,t_q)\in\NN_+^q$, put $\tmax=\max_{j\in[q]}t_j$.
We write $\vt K_r^r = (t_1K_r^r,\dots,t_q K_r^r)$,
where $tK_r^r$ is a matching of size $t$.
Write $\Lambda_\vt=\sum_{i=1}^q(t_i-1)$ for $\vt\in\NN_+^q$.
Generalising the classical theorem of Cockayne and Lorimer~\cite{CL75}
on the multicolour Ramsey number of matchings,
Alon, Frankl and Lov\'asz~\cite{AFL86} proved the following sharp result using topological methods.

\begin{theorem}[Alon--Frankl--Lov\'asz]\label{thm:AFL}
  Fix integers $r\ge 2$ and $q\ge 1$.
  For every $\vt\in\NN_+^q$ and every integer $n\ge(r-1)\tmax+\Lambda_\vt+1$,
  we have $K_n^r\to\vt K_r^r$.
\end{theorem}

\subsection{Main results and applications}\label{sec:intro-results}

Over the past few decades, considerable efforts have been dedicated to establishing sparse random analogues
of classical extremal and Ramsey theorems.
Such \emph{transference} results show that sparse random graphs or hypergraphs inherit
the relevant extremal or Ramsey properties of their complete counterparts.
Examples include sparse random versions of Ramsey's theorem due to R\"odl and Ruci\'nski~\cite{RR95},
and of Tur\'an's theorem due to Conlon and Gowers~\cite{CG16} and Schacht~\cite{Schacht16}.
Gishboliner, Glock, Michaeli and Sgueglia~\cite{GGMS26}
obtained such a transference version of the Alon--Frankl--Lov\'asz theorem for sparse binomial random hypergraphs,
when the target matching sizes are all equal.
Here we seek to isolate the combinatorial properties underlying this transference.
We prove the same asymptotic conclusion for much more general hosts satisfying a weak deterministic condition,
with the random result following as a corollary.

Following~\cite{KM26+}, we say that an $r$-graph $H$ is an \defn{$s$-connector} if,
whenever $A_1,\dots,A_r\subseteq V(H)$ are pairwise disjoint and $|A_i|\ge s$ for every $i\in[r]$,
there exists an edge $e\in E(H)$ such that $|e\cap A_i|=1$ for every $i\in[r]$.
Equivalently, $H$ has \emph{$r$-partite hole number} at most $s-1$ (see~\cites{DT20,BD24}).

\begin{samepage}
  Write $\afl_r(\vt)=(r-1)\tmax+\Lambda_\vt+1$ for the Ramsey number determined by \cref{thm:AFL}.
  Our main result shows that, when $s=o(n)$, the complete $r$-graph in \cref{thm:AFL}
  may be replaced by any $n$-vertex $s$-connector at the cost of $o(n)$ additional vertices.

  \begin{theorem}\label{thm:main}
    For all integers $r\ge 2$ and $q\ge 1$ and every $\eta>0$,
    there exists $\eps>0$ such that the following holds.
    For every $\vt\in\NN_+^q$ and all integers $1\le s\le n$ satisfying $s<\eps n$ and $n\ge\afl_r(\vt)+\eta n$,
    every $n$-vertex $s$-connector $r$-graph $H$ satisfies $H\to\vt K_r^r$.
  \end{theorem}
\end{samepage}

  \Cref{thm:main} resolves a conjecture posed by the authors in~\cite{KM26+}.
  Our earlier proof in the graph case ($r=2$) relies on compression operations based on the Gallai--Edmonds decomposition.
  No analogue of this decomposition is known for general uniform hypergraphs, so we use a different approach (see \cref{sec:intro-methods}).

  As immediate applications,
  \cref{thm:main}
  recovers
  both the defect and transference versions of the Alon--Frankl--Lov\'asz theorem in~\cite{GGMS26}.
  Indeed, an $n$-vertex $r$-graph with at least $(1-\delta)\binom{n}{r}$ edges is an
  $\eps n$-connector, where $\eps=\eps(\delta)\to 0$ as $\delta\to 0$.
  In addition, it is not hard to see that a binomial random $r$-graph
  with edge probability $p\ge Cn^{1-r}$ is,
  \whp{}\footnote{An event holds \defn{with high probability}, or \whp{}, if its probability tends to $1$ as $n\to\infty$.},
  an $\eps n$-connector, where $\eps=\eps(C)\to 0$ as $C\to\infty$.
  More generally, the following theorem extends \cref{thm:main} to random sparsifications of $s$-connectors.

  \begin{theorem}\label{thm:random-thinning}
    For all integers $r\ge 2$ and $q\ge 1$ and every $\eta>0$, there exists $C>0$ such that the following holds.
    For every $\vt\in\NN_+^q$ and all integers $1\le s\le n$ satisfying $n\ge\afl_r(\vt)+\eta n$,
    let $H$ be an $n$-vertex $s$-connector $r$-graph,
    and let $H_p$ be obtained from $H$ by retaining each edge independently with probability $p\in[0,1]$ satisfying $p\ge C(s/n)^{r-1}$.
    With probability at least $1-e^{-n}$, we have $H_p\to\vt K_r^r$ simultaneously for all such $\vt$.
  \end{theorem}

  \Cref{thm:main} also applies to uniformly random regular hypergraphs.
  Fix an integer $r\ge 3$, and let $G$ be chosen uniformly from all simple $d$-regular $r$-graphs on $[n]$,
  where $d=d(n)\to\infty$, $r\mid nd$ and $d\le\binom{n-1}{r-1}$.
  The enumeration and coupling results in~\cites{DFRS13,DFRS17,GIMM23}
  imply that $G$ is, \whp{}, an $s$-connector with $s=o(n)$.
  In particular, for every fixed $q\ge 1$, \whp{} every $q$-colouring of $G$
  contains a monochromatic matching of the asymptotically optimal size $(1/(r+q-1)-o(1))n$.

  A \defn{Steiner triple system} is a $3$-graph in which every pair of vertices lies in exactly one edge.
  DeBiasio and Tait~\cite{DT20}*{Problem~7.4(i)} asked how large a monochromatic matching
  can be guaranteed in a $2$-colouring of a Steiner triple system;
  we give an asymptotically sharp answer for uniformly random systems.
  For $n\equiv 1,3\pmod 6$, a uniformly random Steiner triple system on $n$ vertices is,
  \whp{},
  an $s$-connector with $s=o(n)$~\cites{DT20,FK20}.
  Consequently,
  by \cref{thm:main},
  for every fixed $q\ge 1$, \whp{} every $q$-colouring contains
  a monochromatic matching of the asymptotically optimal size $(1/(q+2)-o(1))n$.

\subsection{Perfect-matching discrepancy}\label{sec:intro-discrepancy}

  The general question in discrepancy theory is whether a ground set can be $q$-coloured so that every set in a given family contains roughly the same number of elements of each colour.
  Here we seek conditions under which every $q$-edge-colouring of an $r$-graph contains a perfect matching with \defn{high discrepancy},
  meaning that significantly more than a $1/q$-proportion of its edges have the same colour.
  Fix $r\ge 3$ and $q\ge 2$, and assume that $r\mid n$.
  By \cref{thm:AFL}, every $q$-colouring of $K_n^r$ has a perfect matching
  with approximately $n/(r+q-1)$ edges of one colour.
  The bound cannot be improved asymptotically for any $n$-vertex $r$-graph,
  as shown by the extremal construction for the complete $r$-graph.

  Write $\bG^{(r)}(n,p)$ for the binomial random $r$-graph on $[n]$,
  in which every $r$-set is an edge independently with probability $p$.
  Johansson, Kahn and Vu~\cite{JKV08+} showed that, for some constant $A=A(r)>0$,
  $\bG^{(r)}(n,p)$ contains a perfect matching \whp{} whenever $p\ge An^{1-r}\log n$.
  Gishboliner, Glock and Sgueglia~\cite{GGS25} asked whether, already at this density,
  \whp{} every $2$-colouring of a random $r$-graph contains a perfect matching
  with approximately $n/(r+1)$ edges of one colour, as in the complete case.
  To our knowledge, this question remains open.
  H\`an, Lang, Marciano, Pavez-Sign\'e, Sanhueza-Matamala, Treglown and Z\'arate-Guer\'en~\cite{HLMPSTZ24+}
  proved that, for $p\ge Cn^{1-r/2}$,
  \whp{} every $q$-colouring of $\bG^{(r)}(n,p)$ contains a perfect matching
  with at least $n/(rq)+\gamma_{r,q}n$ edges of one colour, for some $\gamma_{r,q}>0$.
  They also asked for the correct probability threshold for guaranteeing some positive linear excess over $n/(rq)$ in every $q$-colouring.
  The following proposition gives the asymptotically optimal colour count at $p\ge Cn^{2-r}\log n$,
  a factor $n$ above the perfect-matching existence bound.
  Its proof combines a simple sprinkling argument
  with a large monochromatic matching supplied by \cref{thm:main}.

  \begin{proposition}\label{prop:random-discrepancy}
    For all integers $r\ge 3$ and $q\ge 2$ and every $0<\mu<1$,
    there exists $C>0$ such that the following holds whenever $r\mid n$ and $p\ge Cn^{2-r}\log n$.
    \Whp{}, every $q$-colouring of $G\sim\bG^{(r)}(n,p)$ contains a perfect matching
    with at least $(1-\mu)n/(r+q-1)$ edges of one colour.
  \end{proposition}

  Gishboliner, Glock, Michaeli and Sgueglia~\cite{GGMS26}
  obtain the same colour count already for $p\ge Cn^{1-r}\log n$
  when a $q$-colouring of all $r$-subsets of $[n]$ is fixed before the random hypergraph is sampled.
  Their high-probability conclusion holds for each fixed colouring separately.
  In \cref{prop:random-discrepancy}, an adversary may first observe the sampled hypergraph and only then colour its edges.
  We require, however, the larger density $p\ge Cn^{2-r}\log n$.
  We give the short proof in \cref{sec:main-proofs}.

\subsection{Related work and motivation for the proof}\label{sec:intro-methods}

  There is a substantial literature on transferring Ramsey results for forests and related hypergraph structures
  to sparse binomial random graphs and hypergraphs.
  Such results have been obtained for paths~\cites{Let16,DP17}, trees~\cites{BD17,DP17},
  components and loose paths in hypergraphs~\cite{BDDE19},
  matchings and linear forests~\cite{GKM22a}, and arbitrary forests~\cite{GKM22b}.
  More recent results concern hypergraph matchings~\cite{GGMS26} and forest tilings~\cites{ACFMPY26+,LTWY26+}.
  Liu, Turevskii, Wang and Yan~\cite{LTWY26+} also give a topology-free proof of the Alon--Frankl--Lov\'asz theorem
  for $n$ sufficiently large in terms of the uniformity $r$ and the number of colours $q$,
  extend it to $r$-graphs with sufficiently high minimum degree,
  and give another proof of the random-transference theorem in~\cite{GGMS26}.

  Most of these proofs follow a common approach.\footnote{\label{fn:regularity-exceptions}For approaches not using regularity, see~\cites{ACFMPY26+,GKM22b}.}
  One starts with a pseudorandom host in which the edge density between any $r$ sufficiently large disjoint vertex sets
  is close to a common value $p$, as holds \whp{} for the binomial random model at the relevant densities.
  The multicolour weak sparse hypergraph regularity lemma (see, e.g., in~\cite{GGMS26})
  partitions the vertices into clusters of almost equal size so that almost every $r$-tuple of clusters is regular in every colour.
  The lower-density bound on the host ensures that some colour has density at least $p/(2q)$ between the clusters of each such tuple.
  These tuples form an almost-complete reduced $r$-graph on the clusters, with each edge assigned at least one such colour.

  A defect version of a suitable dense Ramsey theorem then gives a monochromatic structure of asymptotically optimal size in the reduced graph.
  The structure sought there need not be the final object.
  This dense step can require a substantial argument adapted to the structure in question.
  Such defect versions have been obtained, for example, for matchings in graphs~\cite{GKM22a},
  connected matchings~\cite{Let22} and matchings in hypergraphs~\cites{GGMS26,LTWY26+}.
  Finally, regularity and density allow one to lift the reduced structure to the original host.
  For example, a monochromatic connected matching in the reduced graph can be lifted to a long path of the same colour in the original graph~\cite{DP17}.
  For matchings, regularity and density allow us to lift each reduced edge greedily to a matching of the same colour
  covering almost all vertices of its clusters.
  Disjoint reduced edges give disjoint lifted matchings, whose union has asymptotically optimal size.

  Our aim is to weaken these pseudorandomness assumptions substantially,
  retaining only a simple combinatorial condition sufficient for transference.
  In an $n$-vertex $r$-graph, we ask only for a \emph{single} edge crossing any $r$ pairwise disjoint vertex sets
  of size at least $s$, where $s\ll n$ (that is, $s/n\to 0$).
  We refer to such hypergraphs as \defn{sublinear connectors};
  equivalently, these are $r$-graphs with no balanced $r$-partite holes of linear size.
  For Ramsey problems concerning forests, replacing pseudorandomness by this assumption is not new.
  In graphs, Gishboliner, Krivelevich and Michaeli~\cite{GKM22b} showed that sublinear connectors
  guarantee asymptotically as many edges in a monochromatic forest as complete graphs do.
  For two and three colours, Bal and DeBiasio~\cite{BD24} proved the analogous result
  for the order of a largest monochromatic component, and hence for a largest monochromatic tree.

  One might try to imitate the regularity construction in an $s$-connector.
  For simplicity, consider the diagonal case $t_1=\dots=t_q$, with $r,q$ fixed and $s\ll n$.
  Given a $q$-colouring of an $n$-vertex $s$-connector $r$-graph $H$,
  our aim is a monochromatic matching of size $\sim n/(r+q-1)$.
  Partition almost all vertices into $s$-vertex clusters.
  Since $H$ is an $s$-connector, every $r$-tuple of clusters has a crossing edge.
  For each $r$-tuple, we choose one such edge of $H$ and give its colour to the corresponding edge
  of the complete ``reduced'' hypergraph on the clusters.
  By AFL, the reduced hypergraph has a monochromatic matching of about $n/(s(r+q-1))$ edges.
  For each edge of this reduced matching, we take the crossing edge of $H$ that we used to colour it.
  These crossing edges form a monochromatic matching in $H$ of the same size, a factor $s$ short of our target.

  We might instead partition almost all vertices into clusters of size $\ell$, where $s\ll\ell\ll n$.
  Greedily choosing crossing edges between any $r$ of these clusters gives about $\ell$ disjoint edges,
  but these need not have the same colour.
  The complete reduced hypergraph now has about $n/\ell$ vertices, so an almost-perfect matching has about $n/(r\ell)$ edges.
  The crossing matchings for these disjoint $r$-tuples together give about $n/r$ edges.
  The pigeonhole principle then guarantees only about $n/(rq)$ edges of one colour, short of our target when $q\ge 2$.

  We return to $s$-vertex clusters, but now choose the partition more carefully, taking account of the colouring.
  A small vertex cover in one colour would prevent a large matching in that colour.
  Our next aim is to rule out such an obstruction in at least one colour.
  We choose a minimum vertex cover in each colour and group the vertices according to their membership in these $q$ covers.
  Dividing each of the resulting $2^q$ sets into $s$-vertex clusters leaves only $O(s)=o(n)$ vertices unused.
  Each cluster is now either contained in or disjoint from each cover.
  We form the reduced colouring as above.
  A reduced edge of colour $j$ comes from an original edge that meets the colour-$j$ cover,
  so one of its clusters lies entirely inside that cover.
  Thus a colour-$j$ reduced matching of size $k$ gives $k$ disjoint clusters inside the cover, proving that its size is at least $ks$.
  By AFL, some colour has minimum vertex-cover size at least $\sim n/(r+q-1)$.
  A minimum vertex cover can be much larger than a maximum matching, so this does not yet give the desired matching.

  For fractional matchings, however, LP duality states that the maximum fractional matching weight equals the minimum fractional vertex-cover weight.
  We therefore adapt the partition argument to minimum fractional vertex covers,
  using closeness of weights in place of membership in a set.
  We associate each vertex with the vector of its $q$ cover weights in $[0,1]^q$.
  Small boxes in this cube replace the membership classes above: within each box, we divide the vertices into $s$-vertex clusters.
  Since $s\ll n$, we can let the box diameters tend to zero slowly enough that the number of unused vertices remains sublinear.
  All colours must be controlled simultaneously, since the colour supplied by AFL is not known in advance.
  For a crossing edge of colour $j$, the cover weights $a_1,\dots,a_r$ of its vertices satisfy $a_1+\cdots+a_r\ge 1$.
  Near constancy within the clusters then gives total colour-$j$ cover weight at least $s(a_1+\cdots+a_r)\ge s$
  on these clusters, up to the clustering error.
  The total error is therefore sublinear.
  As before, summing over the disjoint clusters of a reduced monochromatic matching gives total cover weight at least $\sim n/(r+q-1)$.
  LP duality then gives
  a monochromatic fractional matching of weight $\sim n/(r+q-1)$.

  To pass to an ordinary matching, we use Kahn's theorem~\cite{Kahn96}.
  It gives a matching of almost the same size if the total fractional weight on edges containing any pair of distinct vertices is $o(1)$.
  Our fractional construction need not have this property.
  For a uniformly random $\ell$-vertex set $U$, a given vertex belongs to $U$ with probability $\ell/n$,
  whereas a given pair belongs to $U$ with probability $\sim\ell^2/n^2$.
  This suggests averaging fractional matchings on small induced subgraphs.
  Each $H[U]$ is still an $s$-connector, so, provided $s\ll\ell$, the fractional result gives
  a monochromatic fractional matching of weight at least $\sim\ell/(r+q-1)$ in every $H[U]$.
  We average these matchings over all $\ell$-sets.
  The average has total weight at least $\sim\ell/(r+q-1)$ and total weight at most $\ell/n$ through any vertex.
  Scaling the average by $n/\ell$ therefore restores the desired total weight while keeping every vertex sum at most $1$.
  Each local fractional matching puts weight at most $1$ through any pair, and a given pair lies in $U$ with probability at most $(\ell/n)^2$.
  After scaling, the total weight through any pair is therefore at most $\ell/n$.
  Taking $\ell\ll n$ makes these pair weights $o(1)$, as required by Kahn's theorem.

  The fractional matchings in different $H[U]$ may have different colours.
  Their average may therefore assign positive weight to edges of several colours.
  To choose a common colour without losing a constant factor, we use concentration.
  By the pigeonhole principle, some colour attains the local bound above on at least a $1/q$ proportion of the $\ell$-sets.
  The maximum fractional matching weight in that colour changes by at most $1$ when one vertex of the set is replaced.
  Azuma's inequality therefore implies concentration within $o(\ell)$ of its mean.
  Since the large value occurs on a positive proportion of the sets, the mean must also be at least $\sim\ell/(r+q-1)$.
  We now average maximum fractional matchings in this fixed colour and scale by $n/\ell$ as above.
  The resulting fractional matching has total weight at least $\sim n/(r+q-1)$ and small pair weights,
  so Kahn's theorem gives an ordinary monochromatic matching of the required size.

\section{Preliminaries}\label{sec:preliminaries}

A \defn{$q$-multi-colouring} of an $r$-graph $H=(V,E)$
is a sequence $\cH=(H_1,\dots,H_q)$
where $H_j=(V,E_j)$ is an $r$-graph for every $j\in[q]$,
and $\bigcup_jE_j=E$.
We call $H$ the \defn{underlying $r$-graph}
of the \defn{coloured $r$-graph} $\cH$.
An edge $e$ is said to have colour $j$ if $e\in E_j$.
Note that we allow an edge to possess multiple colours
and that, for the Ramsey properties considered here, it is equivalent
to consider multi-colourings or usual colourings,
where each edge has exactly one colour.
We will therefore abuse terminology and also
refer to $q$-multi-colourings as \defn{$q$-colourings}.

The \defn{matching number} $\nu(F)$ of a hypergraph $F$
is the size of a maximum matching in $F$.
A \defn{fractional matching} in a hypergraph $F$ is a function $x:E(F)\to[0,\infty)$ such that
\(
  \sum_{e\ni v}x(e)\le 1
\)
for every $v\in V(F)$. Its \defn{weight} is $\sum_{e\in E(F)}x(e)$.
The \defn{fractional matching number} $\nu^*(F)$ is the maximum weight of a fractional matching in $F$.

We will use the following well-known consequence of LP duality to find large fractional matchings.
A \defn{fractional vertex cover} of a hypergraph $F$ is a function $w:V(F)\to[0,\infty)$ such that
\(
  \sum_{v\in e}w(v)\ge 1
\)
for every $e\in E(F)$.
The \defn{weight} of $w$ is the sum of $w(v)$ over all $v\in V(F)$.
The minimum weight of a fractional vertex cover of $F$ is denoted by $\tau^*(F)$.

\begin{lemma}[Fractional matching--cover duality, see, e.g.,~\cite{Berge}]\label{thm:fractional-duality}
  For every hypergraph $F$, we have
  \[
    \nu^*(F)=\tau^*(F).
  \]
\end{lemma}

  We will use the following concentration inequality when passing from fractional to ordinary matchings,
  to show that fractional matching numbers in random induced subhypergraphs are typically close to their means.
  The two-sided form below follows by applying Azuma's inequality (see, e.g.,~\cite{AS}*{Theorem~7.2.1}) to a martingale and its negative.

\begin{lemma}[Two-sided Azuma inequality]\label{thm:Azuma}
  Let $0=X_0,\dots,X_m$ be a martingale satisfying $|X_{i+1}-X_i|\le 1$ for every integer $0\le i<m$.
  For every $\lambda>0$ we have
  \[
    \pr\left(|X_m|>\lambda\sqrt{m}\right)
    <
    2\exp\left(-\frac{\lambda^2}{2}\right).
  \]
\end{lemma}

The starting point of our strategy is the following simple clustering lemma,
applied to the vector obtained from fractional cover weights in each colour,
which will partition the vertices into clusters on most of which the weights are approximately constant.
We write $\diam_\infty$ for diameter in the $\ell_\infty$ metric.

\begin{lemma}[Clustering]\label{lem:blocks}
  For every integer $q\ge 1$ and every $\eta>0$,
  there exists $\eps>0$ such that the following holds.
  Let $V$ be an $n$-element set and let $s\ge 1$ be an integer satisfying $s<\eps n$.
  Let $\rho:V\to[0,1]^q$.
  Put $m=\floor{n/s}$.
  There exists an equipartition $V=V_1\cup\dots\cup V_m$,
  with cluster sizes differing by at most $1$, such that
  \[
    \frac{1}{m}\sum_{i=1}^m\diam_\infty(\rho(V_i))\le\eta.
  \]
\end{lemma}

\begin{proof}
  Partition $[0,1]^q$ into $M=M(\eta,q)$ cubes $W_1,\dots,W_M$ of side length at most $\eta/2$.
  In each coordinate, use intervals closed on the left and open on the right, including $1$ in the last interval.
  Choose an ordering of $V$ in which every element of $\rho^{-1}(W_a)$ precedes every element of $\rho^{-1}(W_b)$ whenever $a<b$.
  Split this order into $m$ consecutive clusters $V_1,\dots,V_m$ whose sizes differ by at most $1$.
  At most $M$ clusters have images meeting more than one cube.

  If a cluster's image lies in one cube, its $\ell_\infty$ diameter is at most $\eta/2$;
  otherwise it is at most $1$.
  Since $s<\eps n$ implies $m>1/\eps-1$,
  choosing $0<\eps<1$ sufficiently small ensures that $M/m\le\eta/2$.
  Therefore
  \[
    \frac{1}{m}\sum_{i=1}^m\diam_\infty(\rho(V_i))\le\frac{\eta}{2}+\frac{M}{m}\le\eta,
  \]
  as required.
\end{proof}

\Needspace{12\baselineskip}
\section{Fractional matchings}\label{sec:fractional}

In this section we prove the following fractional version of \cref{thm:main}.

\begin{samepage}
  \begin{theorem}\label{thm:fractional}
    For all integers $r\ge 2$ and $q\ge 1$ and every $\eta>0$,
    there exists $\eps>0$ such that the following holds.
    Let $1\le s\le n$, let $H$ be an $n$-vertex $s$-connector $r$-graph,
    let $\cH=(H_1,\dots,H_q)$ be a $q$-colouring of $H$,
    and let $\vt\in\NN_+^q$.
    If $s<\eps n$ and $n\ge\afl_r(\vt)+\eta n$,
    then there exists $j\in[q]$ such that $\nu^*(H_j)\ge t_j$.
  \end{theorem}
\end{samepage}

The proof uses the following lemma that clusters vertices according to the
fractional cover weights in each colour.

  \begin{lemma}\label{lem:block-reduction}
    For all integers $r\ge 2$ and $q\ge 1$ and every $\gamma>0$,
    there exists $\eps>0$ such that the following holds.
    Let $H_1,\dots,H_q$ be $r$-graphs on a common $n$-element vertex set $V$,
    and let $s\ge 1$ be an integer satisfying $s<\eps n$.
    Put $m=\floor{n/s}$.
    Then there exists an equipartition $V=V_1\cup\dots\cup V_m$ with the following property.
    For every $j\in[q]$, let $\Gamma_j$ be the $r$-graph on $[m]$ whose edges are the sets
    $S\in\binom{[m]}{r}$ for which $H_j$ has an edge meeting every $V_i$, $i\in S$, once.
    Then, for every $j\in[q]$,
    \[
      \nu^*(H_j)\ge s\nu^*(\Gamma_j)-\gamma n.
    \]
  \end{lemma}

\begin{proof}
    Choose $\eps>0$ as in \cref{lem:blocks} for $q$ and $\gamma$.
    For every $j\in[q]$, choose an optimal fractional vertex cover $w_j$ of $H_j$.
    We may assume that $w_j(v)\le 1$ for every $v\in V$.
    Apply \cref{lem:blocks} to $\rho(v)=(w_1(v),\dots,w_q(v))$
    to obtain an equipartition $V=V_1\cup\dots\cup V_m$ such that
    \[
      \frac{1}{m}\sum_{i=1}^m\diam_\infty(\rho(V_i))\le\gamma.
    \]

    For $i\in[m]$ and $j\in[q]$, put $w_j(i)=\max_{v\in V_i}w_j(v)$.
    For every $v\in V_i$, we have $w_j(i)\le w_j(v)+\diam_\infty(\rho(V_i))$.
    Since $w_j(i)$ does not depend on $v$,
    \[
      w_j(i)
      =\frac{1}{|V_i|}\sum_{v\in V_i}w_j(i)
      \le\frac{1}{|V_i|}\sum_{v\in V_i}\bigl(w_j(v)+\diam_\infty(\rho(V_i))\bigr).
    \]
    Summing over $i$, using $|V_i|\ge s$ and $\sum_vw_j(v)=\nu^*(H_j)$ by \cref{thm:fractional-duality}, gives
    \begin{equation}\label{eq:cover-blocks}
      \begin{aligned}
        \sum_{i=1}^mw_j(i)
        &\le \sum_{i=1}^m\frac{1}{|V_i|}\sum_{v\in V_i}w_j(v)
        +\sum_{i=1}^m\diam_\infty(\rho(V_i))\\
        &\le \frac{\nu^*(H_j)}{s}+\sum_{i=1}^m\diam_\infty(\rho(V_i))
        \le \frac{\nu^*(H_j)}{s}+\gamma m.
      \end{aligned}
    \end{equation}

    Fix $j\in[q]$.
    For every edge $S$ of $\Gamma_j$, choose an edge $e$ of $H_j$ meeting every $V_i$, $i\in S$, once.
    Since $w_j$ is a fractional vertex cover of $H_j$,
    \[
      \sum_{i\in S}w_j(i)\ge\sum_{v\in e}w_j(v)\ge 1.
    \]
    Thus the weights $w_j(i)$ form a fractional vertex cover of $\Gamma_j$.
    By \cref{thm:fractional-duality} and \eqref{eq:cover-blocks},
    \[
      \nu^*(\Gamma_j)
      \le\sum_{i=1}^mw_j(i)
      \le\frac{\nu^*(H_j)}{s}+\gamma m.
    \]
    Rearranging gives $\nu^*(H_j)\ge s\nu^*(\Gamma_j)-\gamma sm\ge s\nu^*(\Gamma_j)-\gamma n$, as required.
\end{proof}

We now deduce \cref{thm:fractional}.  
We use the following properties that are immediate from the definition of $\afl_r$.
Fix integers $r\ge 2$ and $q\ge 1$.
  Write $\one_q=(1,\dots,1)\in\NN^q$.
  For $\vt\in\NN_+^q$ and every integer $b\ge 0$, we have
  \begin{equation}\label{eq:afl-shift}
    \afl_r(\vt+b\one_q)=\afl_r(\vt)+(r+q-1)b.
  \end{equation}
  If $\vt,\vk\in\NN_+^q$ and $0<\theta\le 1$ satisfy $k_j\le\theta t_j$ for every $j\in[q]$, then
  \begin{equation}\label{eq:afl-scaling}
    \afl_r(\vk)\le\theta\afl_r(\vt)+(q-1)(\theta-1)\le\theta\afl_r(\vt).
  \end{equation}

  \begin{proof}[Proof of \cref{thm:fractional}]
    Put $d=r+q-1$ and $\gamma=\eta/(4d)$.
    For these $r$, $q$ and $\gamma$, let $\eps_0>0$ be supplied by \cref{lem:block-reduction}.
    Put $\eps=\min\{\eps_0,\gamma,\eta/(2(d+1))\}$.
    Put $m=\floor{n/s}$ and $b=\ceil{2\gamma n}$.
    Since $b\le 2\gamma n+1$, $s\ge 1$ and $2d\gamma=\eta/2$, we have
    \[
      db+s\le 2d\gamma n+d+s\le\eta n/2+(d+1)s.
    \]
    Moreover, $s<\eps n$ and $\eps\le\eta/(2(d+1))$ give $(d+1)s<\eta n/2$,
    so $db+s<\eta n$.
    Using $n\ge\afl_r(\vt)+\eta n$ and \eqref{eq:afl-shift}, we obtain
    \[
      \afl_r(\vt+b\one_q)=\afl_r(\vt)+db\le n-\eta n+db<n-s\le sm,
    \]
    where the last inequality follows from $m=\floor{n/s}>n/s-1$.
    Apply \cref{lem:block-reduction} to $H_1,\dots,H_q$.
    The resulting graphs form a $q$-colouring $\Gamma=(\Gamma_1,\dots,\Gamma_q)$ of $K_m^r$,
    since $H$ is an $s$-connector.
    Put $k_j=\floor{(t_j+b)/s}$ for every $j\in[q]$ and $\vk=(k_1,\dots,k_q)$.
    Since $b\ge 2\gamma n>2s$, every $k_j$ is positive.
    Applying \eqref{eq:afl-scaling} with $\theta=1/s$ and $\vt+b\one_q$ in place of $\vt$ gives
    \[
      \afl_r(\vk)\le\frac{1}{s}\afl_r(\vt+b\one_q)\le m.
    \]
    By \cref{thm:AFL}, some $\Gamma_j$ has a matching of size $k_j$.
    For this $j$, \cref{lem:block-reduction} and $s<\gamma n$ give
    \[
      \nu^*(H_j)\ge s\nu^*(\Gamma_j)-\gamma n
      \ge sk_j-\gamma n\ge t_j+b-s-\gamma n\ge t_j.\qedhere
    \]
  \end{proof}

\section{Integral matchings}\label{sec:integral}

Now we present the second main ingredient of our proof strategy.
The two key points of the following theorem are
(a) passing from fractional to integral matchings, and
(b) passing from potentially several local choices of $j$ to a single global choice of $j_0$.

\begin{samepage}
  \begin{theorem}[Local-to-global, fractional-to-integral]
    \label{thm:transfer}
    For all integers $r\ge 2$ and $q\ge 1$ and every $\gamma>0$,
    there exist $0<\delta\le 1$ and $\ell_0\in\NN$ such that the following holds.
    Let $n$ and $\ell$ be positive integers satisfying $\ell_0\le\ell<\delta n$,
    let $H_1,\dots,H_q$ be $r$-graphs on a common $n$-element vertex set $V$,
    and let $\vbeta=(\beta_1,\dots,\beta_q)\in[0,1/r]^q$.
    Suppose that, for every $U\in\binom{V}{\ell}$,
    there exists $j\in[q]$ such that $\nu^*(H_j[U])\ge\beta_j\ell$.
    Then there exists $j_0\in[q]$ such that
    \[
      \nu(H_{j_0})\ge(\beta_{j_0}-\gamma)n.
    \]
  \end{theorem}
\end{samepage}

We will address point (a) via the following matching theorem due to Kahn.
For a fractional matching $x$ in a hypergraph $F$, put
\[
  \alpha(x)
  =
  \max_{\substack{u,v\in V(F)\\u\ne v}}
  \sum_{\substack{e\in E(F)\\\{u,v\}\subseteq e}}x(e).
\]

\begin{theorem}[{\cite{Kahn96}*{Theorem~1.2(a)}}]\label{thm:Kahn}
  For every integer $r\ge 2$ and every $\gamma>0$, there exists $0<\delta\le 1$ such that,
  if $F$ is a hypergraph whose edges have size at most $r$
  and $x$ is a fractional matching in $F$ with $\alpha(x)<\delta$,
  then
  \[
    \nu(F)
    \ge
    (1-\gamma)\sum_{e\in E(F)}x(e).
  \]
\end{theorem}

We will address point (b) via the following application of Azuma's inequality,
which provides a single colour $j_0$ with a fractional matching
suitable for applying \cref{thm:Kahn}.

  \begin{lemma}\label{lem:averaging}
    For all integers $r\ge 2$ and $q\ge 1$ and every $\gamma>0$,
    there exists $\ell_0\in\NN$ such that the following holds.
    Let $\ell_0\le\ell<n$, let $H_1,\dots,H_q$ be $r$-graphs on a common $n$-element vertex set $V$,
    and let $\vbeta=(\beta_1,\dots,\beta_q)\in[0,1/r]^q$.
    Suppose that, for every $U\in\binom{V}{\ell}$,
    there exists $j\in[q]$ such that $\nu^*(H_j[U])\ge\beta_j\ell$.
    Then there exist $j_0\in[q]$ and a fractional matching $x$ in $H_{j_0}$ such that
    $\alpha(x)\le(\ell-1)/(n-1)$ and
    \[
      \sum_{e\in E(H_{j_0})}x(e)
      \ge
      (\beta_{j_0}-\gamma)n.
    \]
  \end{lemma}

  The averaging argument in the proof below is similar to that of
  Alon, Frankl, Huang, R\"odl, Ruci\'nski and Sudakov~\cite{AFHRRS12}.

  \begin{proof}
    Choose $\ell_0$ sufficiently large that, for every $\ell\ge\ell_0$,
    $2q\exp(-\gamma^2\ell/2)<1$.
    Let $\bU$ be chosen uniformly from $\binom{V}{\ell}$.
    For $j\in[q]$ and $U\in\binom{V}{\ell}$, put $f_j(U)=\nu^*(H_j[U])$ and $\mu_j=\E f_j(\bU)$.
    We will show that $f_j(\bU)$ is concentrated about its mean $\mu_j$.
    To this end, we apply Azuma's inequality.
    Fix $j\in[q]$ and order the vertices of $\bU$ uniformly at random as $v_1,\dots,v_\ell$.
    For $0\le k\le\ell$, put
    \[
      X_k=\E\bigl[f_j(\bU)\mid v_1,\dots,v_k\bigr]-\mu_j.
    \]
    Then $(X_k)_{k=0}^\ell$ is a martingale with $X_0=0$ and $X_\ell=f_j(\bU)-\mu_j$.

    We need to check that $|X_k-X_{k-1}|\le 1$ for $1\le k\le\ell$.
    Fix $v_1,\dots,v_{k-1}$ and put $W=V\sm\{v_1,\dots,v_{k-1}\}$.
    For $v\in W$, write
    $g(v)=\E\bigl[f_j(\bU)\mid v_1,\dots,v_{k-1},v_k=v\bigr]$.
    We first show that $|g(u)-g(v)|\le 1$ for any $u,v\in W$.
    Pair each possible final set $U$ when $v_k=u$ with $U$ itself if $v\in U$,
    and with $(U\sm\{u\})\cup\{v\}$ otherwise.
    This pairs the possible final sets for the two choices one-to-one; under either choice these sets are uniform.
    Exchanging one vertex of $U$ changes $f_j(U)$ by at most $1$,
    since deleting from a fractional matching all edges through the removed vertex loses weight at most $1$.
    Thus $|g(u)-g(v)|\le 1$.
    Before $v_k$ is exposed, each vertex of $W$ is equally likely to be chosen, so
    \[
      |X_k-X_{k-1}|
      =\left|g(v_k)-\frac{1}{|W|}\sum_{v\in W}g(v)\right|
      \le\frac{1}{|W|}\sum_{v\in W}|g(v_k)-g(v)|
      \le 1.
    \]
    Applying \cref{thm:Azuma} gives
    \[
      \pr\bigl(|f_j(\bU)-\mu_j|>\gamma\ell\bigr)
      <
      2\exp\left(-\frac{\gamma^2\ell}{2}\right).
    \]
    By the union bound over all $j\in[q]$ and the choice of $\ell_0$,
    \[
      \pr\bigl(\exists j\in[q]: |f_j(\bU)-\mu_j|>\gamma\ell\bigr)
      <2q\exp(-\gamma^2\ell/2)<1.
    \]
    Hence there exists an $\ell$-set $U_0$ such that $|f_j(U_0)-\mu_j|\le\gamma\ell$ for every $j\in[q]$.
    By the hypothesis of the lemma, there exists $j_0\in[q]$ such that $f_{j_0}(U_0)\ge\beta_{j_0}\ell$.
    Thus $\mu_{j_0}\ge f_{j_0}(U_0)-\gamma\ell\ge(\beta_{j_0}-\gamma)\ell$.

    For every $\ell$-set $U$, choose an optimal fractional matching $x_U$ in $H_{j_0}[U]$,
    and put $x_U(e)=0$ for every edge $e$ of $H_{j_0}$ not contained in $U$.
    For every $e\in E(H_{j_0})$, put
    \[
      x(e)=\frac{n}{\ell}\E x_{\bU}(e).
    \]
    For every $v\in V$,
    \[
      \sum_{e\ni v}x(e)
      =
      \frac{n}{\ell}\E\sum_{e\ni v}x_{\bU}(e)
      \le
      \frac{n}{\ell}\pr(v\in\bU)
      =
      1,
    \]
    so $x$ is a fractional matching in $H_{j_0}$.
    Moreover,
    \[
      \sum_{e\in E(H_{j_0})}x(e)
      =
      \frac{n}{\ell}\mu_{j_0}
      \ge
      (\beta_{j_0}-\gamma)n.
    \]
    Finally,
    \[
      \alpha(x)
      =\max_{\substack{u,v\in V\\u\ne v}}
      \sum_{e\supseteq\{u,v\}}x(e)
      \le
      \frac{n}{\ell}\max_{\substack{u,v\in V\\u\ne v}}\pr(\{u,v\}\subseteq\bU)
      =
      \frac{\ell-1}{n-1}.\qedhere
    \]
  \end{proof}

\Needspace{6\baselineskip}
  We conclude this section with the proof of its main result.
  
  \begin{proof}[Proof of \cref{thm:transfer}]
    It suffices to consider $0<\gamma\le 1$.
    Apply \cref{thm:Kahn} with $\frac{\gamma}{2}$ in place of $\gamma$,
    and let $0<\delta\le 1$ be the resulting constant.
    Apply \cref{lem:averaging} with $\frac{\gamma}{2}$ in place of $\gamma$,
    and let $\ell_0\in\NN$ be the resulting constant.

    By \cref{lem:averaging}, there exist $j_0\in[q]$ and a fractional matching $x$ in $H_{j_0}$ such that
    \[
      \sum_{e\in E(H_{j_0})}x(e)
      \ge
      \left(\beta_{j_0}-\frac{\gamma}{2}\right)n
    \]
    and $\alpha(x)\le(\ell-1)/(n-1)<\ell/n<\delta$.
    By \cref{thm:Kahn} and $\beta_{j_0}\le 1/r\le 1$, we have
    \[
      \nu(H_{j_0})
      \ge
      \left(1-\frac{\gamma}{2}\right)
      \left(\beta_{j_0}-\frac{\gamma}{2}\right)n
      \ge
      (\beta_{j_0}-\gamma)n.\qedhere
    \]
  \end{proof}

\section{Proofs of the main theorems}\label{sec:main-proofs}

We now assemble the ingredients from the previous sections
to prove the main theorems stated in the introduction.

  \begin{proof}[Proof of \cref{thm:main}]
    It suffices to consider $\eta<1$.
    Put $d=r+q-1$ and $\gamma=\eta/(8d)$.
    For these $r$, $q$ and $\gamma$, let $0<\delta\le 1$ and $\ell_0\in\NN$ be supplied by \cref{thm:transfer},
    and let $\eps_0>0$ be supplied by \cref{lem:block-reduction}.
    Put
    \[
      \eps=\min\left\{\frac{\delta}{2(\ell_0+1)},\frac{\delta\eps_0}{2},\frac{\delta\gamma}{4},\frac{\eta}{2d}\right\}.
    \]

    Let $H=(V,E)$ and $\vt$ satisfy the assumptions of \cref{thm:main},
    and fix a $q$-colouring $\cH=(H_1,\dots,H_q)$ of $H$.
    Put $\ell=\ceil{\delta n}-1$ and $m=\floor{\ell/s}$.
    Since $1\le s<\eps n$, we have $n>1/\eps$.
    By our choice of $\eps$, we have $\delta n/2>\delta/(2\eps)\ge\ell_0+1$,
    so $\ell_0\le\delta n/2\le\ell<\delta n$.
    Our choice of $\eps$ also gives $s<\eps n\le\delta\eps_0 n/2\le\eps_0\ell$.
    Moreover, $s<\eps n\le\delta\gamma n/4\le\gamma\ell/2$,
    and hence $m>\ell/s-1>2/\gamma-1\ge 1/\gamma$.
    Finally, $d\le\eta/(2\eps)<\eta n/2$.

    Put $b=\ceil{4\gamma n}$, so $db\le 4d\gamma n+d=\eta n/2+d<\eta n$.
    Put $\vt^+=\vt+b\one_q$ and $\beta_j=(t_j+b)/n$.
    By~\eqref{eq:afl-shift} and the hypothesis $n\ge\afl_r(\vt)+\eta n$,
    we have $\afl_r(\vt^+)=\afl_r(\vt)+db\le n-\eta n+db<n$.
    For every $j\in[q]$, the definition of $\afl_r$ gives
    \[
      n>\afl_r(\vt^+)=(r-1)\|\vt^+\|_\infty+(t_j+b)+\sum_{i\ne j}(t_i+b-1)\ge r(t_j+b),
    \]
    since $\|\vt^+\|_\infty\ge t_j+b$ and every summand is nonnegative.
    Thus $4\gamma\le b/n\le\beta_j<1/r$ for every $j\in[q]$.
    Put $k_j=\floor{\beta_jm}$ for every $j\in[q]$ and $\vk=(k_1,\dots,k_q)$.
    Since $\beta_jm\ge 4$, every $k_j$ is positive.
    Applying \eqref{eq:afl-scaling} with $\theta=m/n$ and $\vt^+$ in place of $\vt$ gives
    \[
      \afl_r(\vk)\le\frac{m}{n}\afl_r(\vt^+)<m.
    \]

    For every $U\in\binom{V}{\ell}$, apply \cref{lem:block-reduction} to $H_1[U],\dots,H_q[U]$.
    The resulting graphs form a $q$-colouring $\Gamma=(\Gamma_1,\dots,\Gamma_q)$ of $K_m^r$,
    since $H[U]$ is an $s$-connector.
    By \cref{thm:AFL}, there exists $j=j(U)\in[q]$ such that $\Gamma_j$ has a matching of size $k_j$.
    Since $m>1/\gamma$, we have
    $\nu^*(\Gamma_j)\ge k_j\ge\beta_jm-1\ge(\beta_j-\gamma)m$.

    Since $sm\ge\ell-s$, $0\le\beta_j-\gamma\le 1$ and $s<\gamma\ell$,
    \cref{lem:block-reduction} gives
    \[
      \begin{aligned}
        \nu^*(H_j[U])
        &\ge s(\beta_j-\gamma)m-\gamma\ell\\
        &\ge(\beta_j-\gamma)(\ell-s)-\gamma\ell\\
        &=(\beta_j-2\gamma)\ell-(\beta_j-\gamma)s\\
        &\ge(\beta_j-2\gamma)\ell-s>(\beta_j-3\gamma)\ell.
      \end{aligned}
    \]
    Thus \cref{thm:transfer}, with $\beta_j-3\gamma$ in place of $\beta_j$, gives some $j\in[q]$ such that
    \[
      \nu(H_j)\ge(\beta_j-4\gamma)n=t_j+b-4\gamma n\ge t_j.\qedhere
    \]
  \end{proof}

  To prove \cref{thm:random-thinning}, we use the following lemma to verify the connector hypothesis of \cref{thm:main} for the sparsified hypergraph.

  \begin{lemma}\label{lem:sparsification-connector}
    For every integer $r\ge 2$ and every $\alpha>0$, there exists $C>0$ such that the following holds.
    Let $n$ and $s$ be positive integers satisfying $2s\le\alpha n$, let $H$ be an $n$-vertex $s$-connector $r$-graph,
    and let $H_p$ be obtained from $H$ by retaining each edge independently with probability $p\ge C(s/n)^{r-1}$.
    With probability at least $1-e^{-n}$, $H_p$ is a $\ceil{\alpha n}$-connector.
  \end{lemma}

  \begin{proof}
    Put $c=(\alpha/2)^r$ and choose $C>0$ such that $cC\ge\log(r+1)+1$.
    Put $a=\ceil{\alpha n}$.
    Fix pairwise disjoint $a$-sets $A_1,\dots,A_r$.
    We want $H_p$ to contain an edge meeting each $A_i$ once.
    To this end, we first find many such edges in $H$.
    Partition each of $A_1,\dots,A_{r-1}$ into $\floor{a/s}$ disjoint $s$-sets, ignoring the remaining vertices.
    For each $i\in[r-1]$, choose one of these $s$-sets inside $A_i$ and call it $B_i$.
    For each $v\in A_r$, look for an edge consisting of $v$ and one vertex from each $B_i$.
    Such an edge can be absent for at most $s-1$ vertices:
    $s$ of them together with $B_1,\dots,B_{r-1}$ would contradict the connector property.
    Choosing one edge for each of the other vertices gives at least $a-s+1$ distinct edges.
    Each chosen edge has exactly one vertex in each $A_i$.
    For each $i\in[r-1]$, its vertex in $A_i$ lies in exactly one of the disjoint $s$-sets, so this vertex determines $B_i$.
    Hence no edge is counted for two different choices of $B_1,\dots,B_{r-1}$.
    Since $a\ge\alpha n\ge 2s$, we have $a-s+1\ge a/2$ and $\floor{a/s}\ge a/s-1\ge a/(2s)$.
    The number of edges of $H$ meeting every $A_i$ once is therefore at least
    \[
      (a-s+1)\floor{a/s}^{r-1}\ge\frac{a}{2}\left(\frac{a}{2s}\right)^{r-1}
      =\frac{a^r}{2^rs^{r-1}}\ge\frac{cn^r}{s^{r-1}}.
    \]
    Using $p\ge C(s/n)^{r-1}$, the probability that none of these edges survives is at most
    \[
      \exp\left(-\frac{pcn^r}{s^{r-1}}\right)\le\exp(-cCn).
    \]
    There are at most $(r+1)^n$ choices of the disjoint sets $A_1,\dots,A_r$,
    since each vertex lies in one of these $r$ sets or outside their union.
    By the union bound, the probability that $H_p$ is not an $a$-connector is at most
    \[
      (r+1)^n e^{-cCn}=\exp\left((\log(r+1)-cC)n\right)\le e^{-n},
    \]
    where the last inequality uses $cC\ge\log(r+1)+1$.
  \end{proof}

  \begin{proof}[Proof of \cref{thm:random-thinning}]
    Let $0<\eps\le 1$ be the constant supplied by \cref{thm:main} for $r$, $q$ and $\eta$.
    Put $\alpha=\eps/2$, and let $C>0$ be supplied by \cref{lem:sparsification-connector} for $r$ and $\alpha$.
    We may assume that $C\ge(2/\alpha)^{r-1}$.
    Since $C(s/n)^{r-1}\le p\le 1$, we have $s\le C^{-1/(r-1)}n\le\alpha n/2$.
    In particular, $2s\le\alpha n$ and, since $s\ge 1$, $\alpha n\ge 2$.
    Thus $\ceil{\alpha n}\le\alpha n+1<2\alpha n=\eps n$.
    \Cref{lem:sparsification-connector} implies that,
    with probability at least $1-e^{-n}$, $H_p$ is a $\ceil{\alpha n}$-connector.
    On this event, \cref{thm:main} gives $H_p\to\vt K_r^r$
    simultaneously for all $\vt\in\NN_+^q$ satisfying $n\ge\afl_r(\vt)+\eta n$.
  \end{proof}

\begin{samepage}
We prove \cref{prop:random-discrepancy} by using \cref{thm:main} to find a large monochromatic matching and then matching the remaining vertices.
The remaining vertex set depends on the colouring,
so we will show that, \whp{}, every sufficiently large vertex set of size divisible by $r$ has a perfect matching.
We use the following perfect-matching theorem of Johansson, Kahn and Vu~\cite{JKV08+}.

\begin{theorem}[\cite{JKV08+}]\label{thm:random-perfect-matching}
  For every integer $r\ge 2$, there exists $A>0$ such that,
  whenever $r\mid m$ and $p\ge Am^{1-r}\log m$,
  the random $r$-graph $\bG^{(r)}(m,p)$ contains a perfect matching \whp{}.
\end{theorem}
\end{samepage}

We prove the following lemma by a simple sprinkling argument, also used in~\cite{HLMPSTZ24+}.

\begin{lemma}\label{lem:universal-completion}
  For every integer $r\ge 3$ and every $\alpha>0$, there exists $C>0$ such that,
  if $p\ge Cn^{2-r}\log n$ and $G\sim\bG^{(r)}(n,p)$,
  then \whp{} $G[W]$ contains a perfect matching
  for every $W\subseteq[n]$ satisfying $|W|\ge\alpha n$ and $r\mid |W|$.
\end{lemma}

\begin{proof}
  Let $A$ be the constant supplied by \cref{thm:random-perfect-matching},
  and take $C>A\alpha^{1-r}$.
  Choose an integer $f$ with $1\ll f\ll n$, and put $n'=n+f$ and $p_0=p/n'$,
  and let $F_1,\dots,F_{n'}$ be independent copies of $\bG^{(r)}(n,p_0)$ on the same vertex set.
  Fix $W\subseteq[n]$ with $w=|W|\ge\alpha n$ and $r\mid w$.
  Since $n'/n\to 1$, our choice of $C$ gives $p_0\ge A\alpha^{1-r}n^{1-r}\log n\ge Aw^{1-r}\log w$ for all sufficiently large $n$.
  By \cref{thm:random-perfect-matching},
  each $F_i[W]$ has a perfect matching with probability at least $\frac{1}{2}$.
  Hence the probability that none of $F_1[W],\dots,F_{n'}[W]$ has a perfect matching is at most $2^{-n'}$.
  Put $F=F_1\cup\cdots\cup F_{n'}$.
  By the union bound over the at most $2^n$ choices of $W\subseteq[n]$ with $|W|\ge\alpha n$ and $r\mid |W|$,
  the probability that $F[W]$ has no perfect matching for some such $W$ is at most $2^{n-n'}=o(1)$.
  The union $F$ has the distribution of $\bG^{(r)}(n,p')$, where $p'=1-(1-p_0)^{n'}\le n'p_0=p$.
  Coupling with $G$ completes the proof.
\end{proof}

\begin{proof}[Proof of \cref{prop:random-discrepancy}]
  Put $d=r+q-1$ and $t=\ceil{(1-\mu)n/d}$.
  Choose $0<\eps<1/r$ as in \cref{thm:main} for $\eta=\mu/2$,
  and put $s=\ceil{\eps n/2}$.
  The definition of $\afl_r$ gives
  \[
    \afl_r(t\one_q)=(r-1)t+q(t-1)+1=dt-q+1\le dt.
  \]
  Since $t\le(1-\mu)n/d+1$, for all sufficiently large $n$ we have $s<\eps n$ and
  \[
    dt+\frac{\mu n}{2}
    \le
    (1-\mu)n+d+\frac{\mu n}{2}
    =
    (1-\mu/2)n+d
    \le
    n.
  \]
  Thus $n\ge\afl_r(t\one_q)+\mu n/2$.

  Put $\alpha=(q-1)/(2d)$,
  let $C>0$ be supplied by \cref{lem:universal-completion} for $r$ and $\alpha$,
  and let $G\sim\bG^{(r)}(n,p)$ with $p\ge Cn^{2-r}\log n$.
  There are at most $(r+1)^n$ ordered choices of pairwise disjoint sets $A_1,\dots,A_r\subseteq[n]$.
  If each set has size at least $s$, the probability that no edge meets them all is at most $\exp(-ps^r)$.
  Since $(r+1)^n\exp(-ps^r)=o(1)$, a union bound shows that $G$ is an $s$-connector \whp{}.
  By \cref{lem:universal-completion}, \whp{} $G[W]$ also has a perfect matching
  for every $W\subseteq[n]$ with $|W|\ge\alpha n$ and $r\mid |W|$.

  Fix a graph $G$ with both properties and an arbitrary $q$-colouring of $G$.
  By \cref{thm:main}, there exists a monochromatic matching $M_1$ of size $t$;
  if necessary, discard edges so that its size is exactly $t$.
  Put $W=[n]\sm V(M_1)$.
  Then $r\mid |W|$ and, for all sufficiently large $n$,
  \[
    |W|
    =
    n-rt
    \ge
    \frac{q-1}{d}n-r
    \ge
    \alpha n.
  \]
  Hence $G[W]$ has a perfect matching $M_2$.
  The union $M_1\cup M_2$ is a perfect matching and contains at least $t\ge(1-\mu)n/d$ edges of one colour.
\end{proof}

The extra factor $n$ over the perfect-matching threshold comes from the $n'\sim n$
independent copies of $\bG^{(r)}(n,p/n')$ on the same vertex set
in the proof of \cref{lem:universal-completion}:
they make the failure probability small enough to take a union bound over every uncovered vertex set.

\section{Concluding remarks}\label{sec:exact}

  Recall that our earlier graph result~\cite{KM26+} guarantees $H\to\vt K_2^2$
  for every $n$-vertex $s$-connector graph $H$ whenever
  $n\ge\afl_2(\vt)+7(q+1)(s-1)$.
  For $\vt=(t_1,\dots,t_q)\in\NN_+^q$, let $R_s^r(\vt)$ be the least $n$ such that every $q$-colouring of every $n$-vertex $s$-connector $r$-graph
  contains a colour-$j$ matching of size $t_j$ for some $j\in[q]$.
  In this notation, the result of Niu and Wang~\cite{NW26+} states that
  \[
    R_s^2(\vt)
    =
    \Lambda_\vt
    +
    \max\{2s,\tmax+s\}.
  \]

  This suggests the following analogue for $r$-uniform hypergraphs:
  \begin{equation}\label{eq:naive-integral}
    R_s^r(\vt)
    \stackrel{?}{=}
    \Lambda_\vt
    +
    \max\{rs,(r-1)\tmax+s\}.
  \end{equation}

However, this formula is false already for $r=3$.
To see this, let $F$ be the \defn{Fano plane}, the Steiner triple system on seven vertices.

\begin{proposition}\label{prop:fano}
  The \emph{Fano plane} $F$ is a $2$-connector $3$-graph and satisfies $\nu(F)=1$.
\end{proposition}

\begin{proof}
  Any two Fano lines meet, so $\nu(F)=1$.
  To see the connector property, let $A,B,C$ be three pairwise disjoint pairs of Fano points, and let $x$ be the remaining point.
  There are four Fano lines avoiding $x$.
  If none of them has one point in each of $A,B,C$, then each contains both points of one of these pairs.
  Two different lines cannot contain the same pair, since two Fano points determine a unique line.
  The three pairs therefore cannot account for all four lines, a contradiction.
\end{proof}

For $r=3$, $s=2$ and $\vt=(2,1)$, the right-hand side of \eqref{eq:naive-integral} is $7$.
Colouring every edge of $F$ red gives a $\vt K_3^3$-free colouring of a $2$-connector $3$-graph on seven vertices,
so $R_2^3(2,1)\ge 8$.

The right-hand side of \eqref{eq:naive-integral} is nevertheless a lower bound for $R_s^r(\vt)$,
as shown by the following two constructions.

\begin{enumerate}
\item Consider disjoint sets \(B\) of size \(rs-1\)
and \(S_j\) of size \(t_j-1\) for \(j\in[q]\). 
Every \(r\)-set that meets \(S=\bigcup_j S_j\) is an edge,
coloured by some \(j\) for which it meets \(S_j\). 
\item Label the colours so that \(t_1=t_{\max}\). 
Take the complete  \(r\)-graph on the 
disjoint union of sets \(B,S_2,\ldots,S_q\) with
\(|B|=rt_1-1\) and \(|S_j|=t_j-1\) for \(2\le j\le q\).
Add $s-1$ isolated vertices. 
Colour an edge with some \(j\ge2\) if it meets \(S_j\), and otherwise with colour $1$.
\end{enumerate}

Although \eqref{eq:naive-integral} fails, it seems plausible that the additive error above $\afl_r(\vt)$ should remain linear in $s$, as in the $r=2$ case.
The following conjecture proposes such an interpolation between AFL 
and \cref{thm:main} for the full range of $s$.

\begin{conjecture}\label{conj:linear-error}
  For all integers $r\ge 2$ and $q\ge 1$, there exists $C=C(r,q)>0$ such that,
  for every integer $s\ge 1$ and every $\vt\in\NN_+^q$,
  \[
    R_s^r(\vt)\le\afl_r(\vt)+C(s-1).
  \]
\end{conjecture}

\subsection*{Statement on AI use}

OpenAI's Codex was used for editorial assistance, proofreading and literature searches.
The mathematical proof details and exposition are entirely due to the authors.

\bibliography{library}

\begin{bibdiv}
\begin{biblist}

\bib{AFHRRS12}{article}{
      author={Alon, Noga},
      author={Frankl, Peter},
      author={Huang, Hao},
      author={R{\"o}dl, Vojt{\v e}ch},
      author={Ruci{\'n}ski, Andrzej},
      author={Sudakov, Benny},
       title={Large matchings in uniform hypergraphs and the conjectures of
  {E}rd{\H{o}}s and {S}amuels},
        date={2012},
        ISSN={0097-3165},
     journal={Journal of Combinatorial Theory. Series A},
      volume={119},
      number={6},
       pages={1200\ndash 1215},
         url={https://doi.org/10.1016/j.jcta.2012.02.004},
      review={\MR{2915641}},
}

\bib{AFL86}{article}{
      author={Alon, Noga},
      author={Frankl, P\'eter},
      author={Lov\'asz, L\'aszl\'o},
       title={The chromatic number of {K}neser hypergraphs},
        date={1986},
     journal={Transactions of the American Mathematical Society},
      volume={298},
      number={1},
       pages={359\ndash 370},
         url={https://doi.org/10.2307/2000624},
      review={\MR{857448}},
}

\bib{AS}{book}{
      author={Alon, Noga},
      author={Spencer, Joel~H.},
       title={The probabilistic method},
     edition={Fourth},
      series={Wiley Series in Discrete Mathematics and Optimization},
   publisher={John Wiley \& Sons, Inc., Hoboken, NJ},
        date={2016},
        ISBN={978-1-119-06195-3},
      review={\MR{3524748}},
}

\bib{ACFMPY26+}{article}{
      author={Arag{\~a}o, Lucas},
      author={Cheng, Xinbu},
      author={Filipe, Rafael},
      author={Miyazaki, Rafael},
      author={Peng, Danni},
      author={Yan, Zhifei},
       title={{R}amsey properties for tilings in random graphs},
        date={2026-05},
     journal={arXiv e-prints},
      eprint={2605.21471},
         url={https://arxiv.org/abs/2605.21471},
}

\bib{BD17}{article}{
      author={Bal, Deepak},
      author={DeBiasio, Louis},
       title={Partitioning random graphs into monochromatic components},
        date={2017},
     journal={Electronic Journal of Combinatorics},
      volume={24},
      number={1},
       pages={Paper P1.18, 25 pp.},
  url={https://www.combinatorics.org/ojs/index.php/eljc/article/download/v24i1p18/pdf/},
}

\bib{BD24}{article}{
      author={Bal, Deepak},
      author={DeBiasio, Louis},
       title={Large monochromatic components in expansive hypergraphs},
        date={2024},
     journal={Combinatorics, Probability and Computing},
      volume={33},
      number={4},
       pages={467\ndash 483},
         url={https://doi.org/10.1017/S096354832400004X},
}

\bib{BDDE19}{article}{
      author={Bennett, Patrick},
      author={DeBiasio, Louis},
      author={Dudek, Andrzej},
      author={English, Sean},
       title={Large monochromatic components and long monochromatic cycles in
  random hypergraphs},
        date={2019},
     journal={European Journal of Combinatorics},
      volume={76},
       pages={123\ndash 137},
         url={https://doi.org/10.1016/j.ejc.2018.10.001},
}

\bib{Berge}{book}{
      author={Berge, Claude},
       title={Hypergraphs: Combinatorics of finite sets},
      series={North-Holland Mathematical Library},
   publisher={North-Holland, Amsterdam},
        date={1989},
      volume={45},
        ISBN={0-444-87489-5},
}

\bib{CL75}{article}{
      author={Cockayne, Ernest~J.},
      author={Lorimer, Peter~J.},
       title={The {R}amsey number for stripes},
        date={1975},
     journal={Journal of the Australian Mathematical Society},
      volume={19},
      number={2},
       pages={252\ndash 256},
         url={https://doi.org/10.1017/S1446788700029554},
      review={\MR{371733}},
}

\bib{CG16}{article}{
      author={Conlon, David},
      author={Gowers, W.~Timothy},
       title={Combinatorial theorems in sparse random sets},
        date={2016},
     journal={Annals of Mathematics. Second Series},
      volume={184},
      number={2},
       pages={367\ndash 454},
         url={https://doi.org/10.4007/annals.2016.184.2.2},
}

\bib{DT20}{article}{
      author={DeBiasio, Louis},
      author={Tait, Michael},
       title={Large monochromatic components in 3-edge-colored {S}teiner triple
  systems},
        date={2020},
     journal={Journal of Combinatorial Designs},
      volume={28},
      number={6},
       pages={428\ndash 444},
         url={https://doi.org/10.1002/jcd.21707},
}

\bib{DFRS13}{article}{
      author={Dudek, Andrzej},
      author={Frieze, Alan},
      author={Ruci{\'n}ski, Andrzej},
      author={{\v S}ileikis, Matas},
       title={Approximate counting of regular hypergraphs},
        date={2013},
     journal={Information Processing Letters},
      volume={113},
      number={19--21},
       pages={785\ndash 788},
         url={https://doi.org/10.1016/j.ipl.2013.07.018},
}

\bib{DFRS17}{article}{
      author={Dudek, Andrzej},
      author={Frieze, Alan},
      author={Ruci{\'n}ski, Andrzej},
      author={{\v S}ileikis, Matas},
       title={Embedding the {E}rd{\H o}s--{R}{\'e}nyi hypergraph into the
  random regular hypergraph and {H}amiltonicity},
        date={2017},
     journal={Journal of Combinatorial Theory. Series B},
      volume={122},
       pages={719\ndash 740},
         url={https://doi.org/10.1016/j.jctb.2016.09.003},
}

\bib{DP17}{article}{
      author={Dudek, Andrzej},
      author={Pra{\l}at, Pawe{\l}},
       title={On some multicolor {R}amsey properties of random graphs},
        date={2017},
     journal={SIAM Journal on Discrete Mathematics},
      volume={31},
      number={3},
       pages={2079\ndash 2092},
         url={https://doi.org/10.1137/16M1069717},
}

\bib{FK20}{article}{
      author={Ferber, Asaf},
      author={Kwan, Matthew},
       title={Almost all {S}teiner triple systems are almost resolvable},
        date={2020},
     journal={Forum of Mathematics, Sigma},
      volume={8},
       pages={e39, 24},
         url={https://doi.org/10.1017/fms.2020.29},
}

\bib{GGMS26}{article}{
      author={Gishboliner, Lior},
      author={Glock, Stefan},
      author={Michaeli, Peleg},
      author={Sgueglia, Amedeo},
       title={Defect and transference versions of the
  {A}lon--{F}rankl--{L}ov{\'a}sz theorem},
        date={2026},
     journal={Combinatorics, Probability and Computing},
      volume={35},
       pages={542\ndash 553},
         url={https://doi.org/10.1017/S0963548326100431},
}

\bib{GGS25}{article}{
      author={Gishboliner, Lior},
      author={Glock, Stefan},
      author={Sgueglia, Amedeo},
       title={Tight {H}amilton cycles with high discrepancy},
        date={2025},
     journal={Combinatorics, Probability and Computing},
      volume={34},
      number={4},
       pages={565\ndash 584},
         url={https://doi.org/10.1017/S0963548325000057},
}

\bib{GKM22a}{article}{
      author={Gishboliner, Lior},
      author={Krivelevich, Michael},
      author={Michaeli, Peleg},
       title={Color-biased {H}amilton cycles in random graphs},
        date={2022},
     journal={Random Structures \& Algorithms},
      volume={60},
      number={3},
       pages={289\ndash 307},
         url={https://doi.org/10.1002/rsa.21043},
}

\bib{GKM22b}{article}{
      author={Gishboliner, Lior},
      author={Krivelevich, Michael},
      author={Michaeli, Peleg},
       title={Discrepancies of spanning trees and {H}amilton cycles},
        date={2022},
     journal={Journal of Combinatorial Theory. Series B},
      volume={154},
       pages={262\ndash 291},
         url={https://doi.org/10.1016/j.jctb.2022.01.003},
}

\bib{GIMM23}{article}{
      author={Greenhill, Catherine},
      author={Isaev, Mikhail},
      author={Makai, Tam{\'a}s},
      author={McKay, Brendan~D.},
       title={Degree sequences of sufficiently dense random uniform
  hypergraphs},
        date={2023},
     journal={Combinatorics, Probability and Computing},
      volume={32},
      number={2},
       pages={183\ndash 224},
         url={https://doi.org/10.1017/S0963548322000190},
}

\bib{HLMPSTZ24+}{article}{
      author={H\`an, Hi\d{\^e}p},
      author={Lang, Richard},
      author={Marciano, Jo\~ao~Pedro},
      author={Pavez-Sign\'e, Mat\'ias},
      author={Sanhueza-Matamala, Nicol\'as},
      author={Treglown, Andrew},
      author={Z\'arate-Guer\'en, Camila},
       title={Colour-bias perfect matchings in hypergraphs},
        date={2025},
     journal={SIAM Journal on Discrete Mathematics},
      volume={39},
      number={3},
       pages={1894\ndash 1916},
         url={https://doi.org/10.1137/24M170483X},
}

\bib{JKV08+}{article}{
      author={Johansson, Anders},
      author={Kahn, Jeff},
      author={Vu, Van},
       title={Factors in random graphs},
        date={2008},
     journal={Random Structures \& Algorithms},
      volume={33},
      number={1},
       pages={1\ndash 28},
         url={https://doi.org/10.1002/rsa.20224},
}

\bib{Kahn96}{article}{
      author={Kahn, Jeff},
       title={A linear programming perspective on the
  {F}rankl--{R}{\"o}dl--{P}ippenger theorem},
        date={1996},
     journal={Random Structures \& Algorithms},
      volume={8},
      number={2},
       pages={149\ndash 157},
      review={\MR{1607104}},
}

\bib{KM26+}{article}{
      author={{Keevash}, Peter},
      author={{Michaeli}, Peleg},
       title={{A very robust Ramsey theorem for matchings}},
        date={2026-03},
     journal={arXiv e-prints},
       pages={arXiv:2603.03139},
      eprint={2603.03139},
}

\bib{Let16}{article}{
      author={Letzter, Shoham},
       title={Path {R}amsey number for random graphs},
        date={2016},
     journal={Combinatorics, Probability and Computing},
      volume={25},
      number={4},
       pages={612\ndash 622},
         url={https://doi.org/10.1017/S0963548315000279},
}

\bib{Let22}{article}{
      author={Letzter, Shoham},
       title={An improvement on {{\L}uczak}'s connected matchings method},
        date={2022},
     journal={Bulletin of the London Mathematical Society},
      volume={54},
      number={2},
       pages={609\ndash 623},
         url={https://doi.org/10.1112/blms.12587},
}

\bib{LTWY26+}{article}{
      author={Liu, Hong},
      author={Turevskii, Maksim},
      author={Wang, Lanchao},
      author={Yan, Zhifei},
       title={Structural reductions for monochromatic matchings and {R}amsey
  tilings},
        date={2026-06},
     journal={arXiv e-prints},
       pages={arXiv:2606.24863},
      eprint={2606.24863},
         url={https://arxiv.org/abs/2606.24863},
}

\bib{NW26+}{article}{
      author={{Niu}, Mengyuan},
      author={{Wang}, Lanchao},
       title={{An exact robust Ramsey theorem for matchings}},
        date={2026-06},
     journal={arXiv e-prints},
       pages={arXiv:2606.19851},
      eprint={2606.19851},
}

\bib{Ramsey30}{article}{
      author={Ramsey, Frank~P.},
       title={On a problem of formal logic},
        date={1930},
     journal={Proceedings of the London Mathematical Society. Second Series},
      volume={30},
      number={1},
       pages={264\ndash 286},
         url={https://doi.org/10.1112/plms/s2-30.1.264},
}

\bib{RR95}{article}{
      author={R{\"o}dl, Vojt{\v e}ch},
      author={Ruci{\'n}ski, Andrzej},
       title={Threshold functions for {R}amsey properties},
        date={1995},
     journal={Journal of the American Mathematical Society},
      volume={8},
      number={4},
       pages={917\ndash 942},
         url={https://doi.org/10.1090/S0894-0347-1995-1276825-6},
}

\bib{Schacht16}{article}{
      author={Schacht, Mathias},
       title={Extremal results for random discrete structures},
        date={2016},
     journal={Annals of Mathematics. Second Series},
      volume={184},
      number={2},
       pages={333\ndash 365},
         url={https://doi.org/10.4007/annals.2016.184.2.1},
}

\end{biblist}
\end{bibdiv}

\end{document}